\documentclass{article}

\usepackage{ucs}
\usepackage[utf8]{inputenc}
\usepackage[greek,english]{babel}
\newcommand{\en}{\selectlanguage{english}}

\begin{document}

\en

\title{Diffeomorphisms of the closed unit disc converging to the identity\footnote{Dedicated to the memory of my grandparents Nikolaos and Alexandra, and Konstantinos and Eleni.}}

\author{Nikolaos E. Sofronidis\footnote{$A \Sigma MA:$ 130/2543/94}}

\date{\footnotesize Department of Economics, University of Ioannina, Ioannina 45110, Greece.
(nsofron@otenet.gr, nsofron@cc.uoi.gr)}

\maketitle

\begin{abstract}
If $\mathcal{G}$ is the group (under composition) of diffeomorphisms $f : {\overline{D}}(0;1) \rightarrow {\overline{D}}(0;1)$ of the closed unit disc ${\overline{D}}(0;1)$ which are the identity map $id : {\overline{D}}(0;1) \rightarrow {\overline{D}}(0;1)$ on the closed unit circle and satisfy the condition $det(J(f)) > 0$, where $J(f)$ is the Jacobian matrix of $f$ or (equivalently) the Fr\'{e}chet derivative of $f$, then $\mathcal{G}$ equipped with the metric $d_{\mathcal{G}}(f,g) = \Vert f-g \Vert _{\infty } + \Vert J(f) - J(g) \Vert _{\infty }$, where $f$, $g$ range over $\mathcal{G}$, is a metric space in which $d_{\mathcal{G}} \left( f_{t} , id \right) \rightarrow 0$ as $t \rightarrow 1^{+}$, where $f_{t}(z) = \frac{ tz }{ 1 + (t-1) \vert z \vert }$, whenever $z \in {\overline{D}}(0;1)$ and $t \geq 1$.
\end{abstract}

\section*{\footnotesize{{\bf Mathematics Subject Classification:} 03E65, 22A05, 22A10, 54E35, 54E40, 54H11.}}

\section{Introduction}

If $\mathcal{G}$ is the group (under composition) of diffeomorphisms $f : {\overline{D}}(0;1) \rightarrow {\overline{D}}(0;1)$ of the closed unit disc ${\overline{D}}(0;1)$ which are the identity map $id : {\overline{D}}(0;1) \rightarrow {\overline{D}}(0;1)$ on the closed unit circle and satisfy the condition $det(J(f)) > 0$, where $J(f)$ is the Jacobian matrix of $f$ or (equivalently) the Fr\'{e}chet derivative of $f$, then $\mathcal{G}$ equipped with the metric $$d_{\mathcal{G}}(f,g) = \Vert f-g \Vert _{\infty } + \Vert J(f) - J(g) \Vert _{\infty },$$ where $f$, $g$ range over $\mathcal{G}$, is a metric space in ' {\it ZF - Axiom of Foundation + Axiom of Countable Choice} '. Our purpose in this article is to prove in ' {\it ZF - Axiom of Foundation + Axiom of Countable Choice} ' the following result.
\\ \rm \\
{\bf 1.1. Theorem.} If $f_{t}(z) = \frac{ tz }{ 1 + (t-1) \vert z \vert }$, whenever $z \in {\overline{D}}(0;1)$ and $t \geq 1$, then $f_{t} \in \mathcal{G}$ and $d_{\mathcal{G}} \left( f_{t} , id \right) \rightarrow 0$ as $t \rightarrow 1^{+}$.

\section{The maximums of the quadratic functions}

Let $a$, $b$, $c$ be arbitrary real numbers such that $b \neq 0$ and either $a \neq 0$ or $c \neq 0$, while $$F : {\bf R}^{2} \ni (u,v) \mapsto (au+bv)^{2} + (bu+cv)^{2} \in {\bf R}.$$ We want to find the maximum of $F$ on the closed unit disc ${\overline{D}}(0;1)$. By virtue of 42.12 on pages 406-407 and 42.13 on page 407 of [1], if $$G : {\bf R}^{2} \ni (u,v) \mapsto u^{2} + v^{2} - 1 \in {\bf R},$$ we are looking for solutions in the unknowns $u$, $v$, $\lambda $ of the following system of equations
\[ \left\{ \begin{array}{lllll} \frac{ \partial F}{\partial u} = \lambda \frac{ \partial G}{\partial u} \\ \\ \frac{ \partial F}{\partial v} = \lambda \frac{ \partial G}{\partial v} \\ \\ G(u,v) = 0 \end{array} \right. \]
Since $$\frac{ \partial F}{\partial u} = 2(au+bv)a + 2(bu+cv)b = 2 \left( \left( a^{2} + b^{2} \right) u + b(a+c)v \right) $$ and $$\frac{ \partial F}{\partial v} = 2(au+bv)b + 2(bu+cv)c = 2 \left( b(a+c)u + \left( b^{2} + c^{2} \right) v \right) $$ the system above takes the form
\[ \left\{ \begin{array}{lllll} \left( a^{2} + b^{2} - \lambda \right) u + (a+c)bv = 0 \\ \\ (a+c)bu + \left( b^{2} + c^{2} - \lambda \right) v = 0 \\ \\ u^{2} + v^{2} = 1 \end{array} \right. \]
The condition $u^{2} + v^{2} = 1$ implies that $(u,v) \neq (0,0)$ and consequently the following determinant must vanish
\[ \left\vert \begin{array}{lll} a^{2} + b^{2} - \lambda & (a+c)b \\ \\ (a+c)b & b^{2} + c^{2} - \lambda \end{array} \right\vert = 0 \]
or (equivalently) $${\lambda }^{2} - \left( a^{2} + 2b^{2} + c^{2} \right) \lambda + \left( b^{2} - ac \right) ^{2} = 0.$$ So $$\lambda = \frac{ a^{2} + 2b^{2} + c^{2} \pm \sqrt{ \left( a^{2} - c^{2} \right) ^{2} + 4b^{2}(a+c)^{2} }}{2}$$ and we distinguish the following three cases:

\subsection{$a^{2}>c^{2}$}

If $$\lambda = \frac{ a^{2} + 2b^{2} + c^{2} + \sqrt{ \left( a^{2} - c^{2} \right) ^{2} + 4b^{2}(a+c)^{2} }}{2},$$ then the above system of three equations implies that $$u = \pm \frac{ \lambda - b^{2} - c^{2} }{ \sqrt{ \left( b^{2} + c^{2} - \lambda \right) ^{2} + (a+c)^{2}b^{2} } }$$ and $$v = \pm \frac{ (a+c)b }{ \sqrt{ \left( b^{2} + c^{2} - \lambda \right) ^{2} + (a+c)^{2}b^{2} } }$$ hence
\begin{enumerate}
\item[ ]
$F(u,v) = \frac{ \left( a^{2} - c^{2} + \sqrt{ \left( a^{2} - c^{2} \right) ^{2} + 4b^{2}(a+c)^{2} } \right) ^{2} \left( a^{2} + b^{2} \right) }{ \left( a^{2} - c^{2} + \sqrt{ \left( a^{2} - c^{2} \right) ^{2} + 4b^{2}(a+c)^{2} } \right) ^{2} + 4(a+c)^{2}b^{2} }$
\item[ ]
$+ \frac{ 4 \left( a^{2} + b^{2} + \sqrt{ \left( a^{2} - c^{2} \right) ^{2} + 4b^{2}(a+c)^{2} } \right) \cdot b^{2}(a+c)^{2} }{ \left( a^{2} - c^{2} + \sqrt{ \left( a^{2} - c^{2} \right) ^{2} + 4b^{2}(a+c)^{2} } \right) ^{2} + 4(a+c)^{2}b^{2} }$,
\end{enumerate}
while if $$\lambda = \frac{ a^{2} + 2b^{2} + c^{2} - \sqrt{ \left( a^{2} - c^{2} \right) ^{2} + 4b^{2}(a+c)^{2} }}{2},$$ then the above system of three equations implies that $$v = \pm \frac{ a^{2} + b^{2} - \lambda }{ \sqrt{ \left( a^{2} + b^{2} - \lambda \right) ^{2} + (a+c)^{2}b^{2} } }$$ and $$u = \mp \frac{ (a+c)b }{ \sqrt{ \left( a^{2} + b^{2} - \lambda \right) ^{2} + (a+c)^{2}b^{2} } }$$ hence
\begin{enumerate}
\item[ ]
$F(u,v) = \frac{ \left( a^{2} - c^{2} + \sqrt{ \left( a^{2} - c^{2} \right) ^{2} + 4b^{2}(a+c)^{2} } \right) ^{2} \left( b^{2} + c^{2} \right) }{ \left( a^{2} - c^{2} + \sqrt{ \left( a^{2} - c^{2} \right) ^{2} + 4b^{2}(a+c)^{2} } \right) ^{2} + 4(a+c)^{2}b^{2} }$
\item[ ]
$+ \frac{ 4 \left( b^{2} + c^{2} - \sqrt{ \left( a^{2} - c^{2} \right) ^{2} + 4b^{2}(a+c)^{2} } \right) b^{2}(a+c)^{2} }{ \left( a^{2} - c^{2} + \sqrt{ \left( a^{2} - c^{2} \right) ^{2} + 4b^{2}(a+c)^{2} } \right) ^{2} + 4(a+c)^{2}b^{2} }$.
\end{enumerate}
Therefore, if we set $$A = \left( a^{2} - c^{2} + \sqrt{ \left( a^{2} - c^{2} \right) ^{2} + 4b^{2}(a+c)^{2} } \right) ^{2} > 0$$ and $$B = 4(a+c)^{2}b^{2} > 0,$$ then the maximum of $F$ on ${\overline{D}}(0;1)$ is one of the numbers
\begin{enumerate}
\item[ ]
$\frac{A}{A+B} \left( a^{2} + b^{2} \right) + \frac{B}{A+B} \left( a^{2} + b^{2} + \sqrt{ \left( a^{2} - c^{2} \right) ^{2} + 4b^{2}(a+c)^{2} } \right) $,
\item[ ]
$\frac{A}{A+B} \left( b^{2} + c^{2} \right) + \frac{B}{A+B} \left( b^{2} + c^{2} - \sqrt{ \left( a^{2} - c^{2} \right) ^{2} + 4b^{2}(a+c)^{2} } \right) $.
\end{enumerate}

\subsection{$a^{2}<c^{2}$}

By symmetry, the same conclusion is true as in the previous subsection with $a$, $c$ interchanged.

\subsection{$a^{2}=c^{2}$}

If $a=c$, then $F(u,v) = \left( a^{2} + b^{2} \right) \left( u^{2} + v^{2} \right) + 4abuv \leq a^{2} + b^{2} + 4abuv$, whenever $(u,v) \in {\overline{D}}(0;1)$, so the maximum of $F$ is one of the numbers $a^{2} + b^{2} \pm 2ab$, while if $a=-c$, then $F(u,v) = \left( a^{2} + b^{2} \right) \left( u^{2} + v^{2} \right) \leq a^{2} + b^{2}$, whenever $(u,v) \in {\overline{D}}(0;1)$, so the maximum of $F$ is the number $a^{2} + b^{2}$.

\section{The computation of the Jacobians}

Let $$f^{(t)}(x,y) = \left( \frac{tx}{ 1 + (t-1) \sqrt{ x^{2} + y^{2} } } , \frac{ty}{ 1 + (t-1) \sqrt{ x^{2} + y^{2} } } \right) ,$$ where $(x,y) \in {\bf R}^{2}$ and $x^{2} + y^{2} \leq 1$, while $t>1$. If $(x,y) \neq (0,0)$ and $$f^{(t)} = \left( f_{1}^{(t)} , f_{2}^{(t)} \right) ,$$ then a straightforward computation shows that $$\left( f_{1}^{(t)} \right) _{x}' = \frac{ t \sqrt{ x^{2} + y^{2} } + t(t-1)y^{2} }{ \sqrt{ x^{2} + y^{2} } \left( 1 + (t-1) \sqrt{ x^{2} + y^{2} } \right) ^{2} } ,$$ $$\left( f_{1}^{(t)} \right) _{y}' = \left( f_{2}^{(t)} \right) _{x}' = \frac{ t(1-t)xy }{ \sqrt{ x^{2} + y^{2} } \left( 1 + (t-1) \sqrt{ x^{2} + y^{2} } \right) ^{2} } $$ and $$\left( f_{2}^{(t)} \right) _{y}' = \frac{ t \sqrt{ x^{2} + y^{2} } + t(t-1)x^{2} }{ \sqrt{ x^{2} + y^{2} } \left( 1 + (t-1) \sqrt{ x^{2} + y^{2} } \right) ^{2} } .$$ Therefore, $$J \left( f^{(t)} \right) = \frac{1}{ \sqrt{ x^{2} + y^{2} } \left( 1 + (t-1) \sqrt{ x^{2} + y^{2} } \right) ^{2} } \cdot $$ \[ \cdot \left[ \begin{array}{lll} t \sqrt{ x^{2} + y^{2} } + t(t-1)y^{2}  &  t(1-t)xy  \\  \\ t(1-t)xy  &  t \sqrt{ x^{2} + y^{2} } + t(t-1)x^{2} \end{array} \right] \] and consequently for any $(u,v) \in {\bf R}^{2}$, we have that $$\left( J \left( f^{(t)} \right) (x,y) \right) \left[ \begin{array}{lll} u \\ \\ v \end{array} \right] = \frac{1}{ \sqrt{ x^{2} + y^{2} } \left( 1 + (t-1) \sqrt{ x^{2} + y^{2} } \right) ^{2} } \cdot $$  \[ \cdot \left[ \begin{array}{lll} \left( t \sqrt{ x^{2} + y^{2} } + t(t-1)y^{2} \right) u &  \left( t(1-t)xy \right) v \\  \\ \left( t(1-t)xy \right) u & \left( t \sqrt{ x^{2} + y^{2} } + t(t-1)x^{2} \right) v \end{array} \right] \] which implies that $$\left\Vert \left( J \left( f^{(t)} \right) (x,y) \right) \left[ \begin{array}{lll} u \\ \\ v \end{array} \right] \right\Vert ^{2} = \frac{ (au+bv)^{2} + (bu+cv)^{2} }{ \left( x^{2} + y^{2} \right) \left( 1 + (t-1) \sqrt{ x^{2} + y^{2} } \right) ^{4} } ,$$ where $a = t \sqrt{ x^{2} + y^{2} } + t(t-1)y^{2}$, $b = t(1-t)xy$ and $c = t \sqrt{ x^{2} + y^{2} } + t(t-1)x^{2}$, and hence $a-c = t(t-1) \left( y^{2} - x^{2} \right) $, $a+c = 2t \sqrt{ x^{2} + y^{2} } + t(t-1) \left( x^{2} + y^{2} \right) $, $a^{2} - c^{2} = t^{2}(t-1) \left( y^{2} - x^{2} \right) \left( 2 \sqrt{ x^{2} + y^{2} } + (t-1) \left( x^{2} + y^{2} \right) \right) $, $(a-c)^{2} + 4b^{2} = t^{2}(t-1)^{2} \left( y^{2} - x^{2} \right) ^{2} + 4t^{2}(1-t)^{2}x^{2}y^{2} = t^{2}(t-1)^{2} \left( x^{2} + y^{2} \right) ^{2}$. Thus, we obtain that $B = 4(a+c)^{2}b^{2} = 4t^{2} \left( 2 \sqrt{ x^{2} + y^{2} } + (t-1) \left( x^{2} + y^{2} \right) \right) ^{2} \cdot t^{2}(1-t)^{2}x^{2}y^{2} = 4t^{4}(t-1)^{2}x^{2}y^{2} \cdot \left( 2 \sqrt{ x^{2} + y^{2} } + (t-1) \left( x^{2} + y^{2} \right) \right) ^{2}$, while
\begin{enumerate}
\item[ ]
$A = \left( a^{2} - c^{2} + \sqrt{ \left( a^{2} - c^{2} \right) ^{2} + 4b^{2}(a+c)^{2} } \right) ^{2}$
\item[ ]
$= \left( a^{2} - c^{2} + \sqrt{ (a+c)^{2} \left( (a-c)^{2} + 4b^{2} \right) } \right) ^{2}$
\item[ ]
$= \left( a^{2} - c^{2} + \vert a+c \vert \cdot \sqrt{ (a-c)^{2} + 4b^{2} } \right) ^{2}$,
\end{enumerate}
where $\vert a+c \vert = a+c$, since $x^{2}+y^{2} \leq 1$ and hence we have that $$a+c = t \sqrt{ x^{2} + y^{2} } \left( 2 + (t-1) \sqrt{ x^{2} + y^{2} } \right) \geq 0,$$ which implies that
\begin{enumerate}
\item[ ]
$A = \left( a^{2} - c^{2} + (a+c) \sqrt{ (a-c)^{2} + 4b^{2} } \right) ^{2} $
\item[ ]
$= \left( t^{2}(t-1) \left( y^{2} - x^{2} \right) \left( 2 \sqrt{ x^{2} + y^{2} } + (t-1) \left( x^{2} + y^{2} \right) \right) \right.$
\item[ ]
$\left. + \left( 2t \sqrt{ x^{2} + y^{2} } + t(t-1) \left( x^{2} + y^{2} \right) \right) \cdot t \vert t-1 \vert \left( x^{2} + y^{2} \right) \right) ^{2}$
\item[ ]
$= \left( t^{2} \left( (t-1) \left( y^{2} - x^{2} \right) + \vert t-1 \vert \left( x^{2} + y^{2} \right) \right) \right.$
\item[ ]
$\left. \cdot \left( 2 \sqrt{ x^{2} + y^{2} } + (t-1) \left( x^{2} + y^{2} \right) \right) \right) ^{2}$
\item[ ]
$= t^{4} (t-1)^{2} \left( x^{2} + y^{2} + sgn(t-1) \cdot \left( y^{2} - x^{2} \right) \right) ^{2}$
\item[ ]
$\cdot \left( 2 \sqrt{ x^{2} + y^{2} } + (t-1) \left( x^{2} + y^{2} \right) \right) ^{2}$
\item[ ]
$= t^{4} (t-1)^{2} \cdot \left( \left( 1 - sgn(t-1) \right) x^{2} + \left( 1 + sgn(t-1) \right) y^{2} \right) ^{2}$
\item[ ]
$\cdot \left( 2 \sqrt{ x^{2} + y^{2} } + (t-1) \left( x^{2} + y^{2} \right) \right) ^{2}$
\end{enumerate}
and
\begin{enumerate}
\item[ ]
$A + B = t^{4}(t-1)^{2} \cdot \left( 2 \sqrt{ x^{2} + y^{2} } + (t-1) \left( x^{2} + y^{2} \right) \right) ^{2}$
\item[ ]
$\cdot \left( \left( \left( 1 - sgn(t-1) \right) x^{2} + \left( 1 + sgn(t-1) \right) y^{2} \right) ^{2} + 4x^{2}y^{2} \right) $.
\end{enumerate}
Therefore, we have that $$\frac{A}{A+B} = \frac{ \left( \left( 1 - sgn(t-1) \right) x^{2} + \left( 1 + sgn(t-1) \right) y^{2} \right) ^{2} }{ \left( \left( 1 - sgn(t-1) \right) x^{2} + \left( 1 + sgn(t-1) \right) y^{2} \right) ^{2} + 4x^{2}y^{2} }$$ and $$\frac{B}{A+B} = \frac{ 4x^{2}y^{2} }{ \left( \left( 1 - sgn(t-1) \right) x^{2} + \left( 1 + sgn(t-1) \right) y^{2} \right) ^{2} + 4x^{2}y^{2} }.$$ It is not difficult to see that $\frac{A}{A+B} = \frac{ y^{2} }{ x^{2} + y^{2} }$ and $\frac{B}{A+B} = \frac{ x^{2} }{ x^{2} + y^{2} }$. Moreover, both $a$, $c$ are positive and consequently the case $a=-c$ can not occur. Since $c$ is simply $a$ with $x$, $y$ interchanged, we need only treat the case $a>c$, among the cases $a>c$ and $a<c$, for if $a<c$, then the result will be as in the case $a>c$ with $x$, $y$ interchanged. We will treat the case $a=c$ at the end. So let us assume that $a>c$ or (equivalently) that $y^{2}>x^{2}$. Then
\begin{enumerate}
\item[ ]
$\sqrt{ \left( a^{2} - c^{2} \right) ^{2} + 4b^{2}(a+c)^{2} } = (a+c) \sqrt{ (a-c)^{2} + 4b^{2} }$
\item[ ]
$= \left( 2t \sqrt{ x^{2} + y^{2} } + t(t-1) \left( x^{2} + y^{2} \right) \right) \cdot t(t-1) \left( x^{2} + y^{2} \right) $
\item[ ]
$= t^{2}(t-1) \left( x^{2} + y^{2} \right) \left( 2 \sqrt{ x^{2} + y^{2} } + (t-1) \left( x^{2} + y^{2} \right) \right) $,
\end{enumerate}
\begin{enumerate}
\item[ ]
$a^{2} + b^{2} = \left( t \sqrt{ x^{2} + y^{2} } + t(t-1)y^{2} \right) ^{2} + t^{2}(1-t)^{2}x^{2}y^{2}$
\item[ ]
$= t^{2} \left( x^{2} + y^{2} \right) + 2t^{2}(t-1)y^{2} \sqrt{ x^{2} + y^{2} } + t^{2}(t-1)^{2}y^{4} + t^{2}(1-t)^{2}x^{2}y^{2}$
\item[ ]
$= t^{2}(t-1)^{2} \left( y^{4} + x^{2}y^{2} \right) + 2t^{2}(t-1)y^{2} \sqrt{ x^{2} + y^{2} } + t^{2} \left( x^{2} + y^{2} \right) $,
\end{enumerate}
\begin{enumerate}
\item[ ]
$b^{2} + c^{2} = t^{2}(1-t)^{2}x^{2}y^{2} + \left( t \sqrt{ x^{2} + y^{2} } + t(t-1)x^{2} \right) ^{2}$
\item[ ]
$= t^{2}(1-t)^{2}x^{2}y^{2} + t^{2} \left( x^{2} + y^{2} \right) + 2t^{2}(t-1)x^{2} \sqrt{ x^{2} + y^{2} } + t^{2}(t-1)^{2}x^{4}$
\item[ ]
$= t^{2}(t-1)^{2} \left( x^{4} + x^{2}y^{2} \right) + 2t^{2}(t-1)x^{2} \sqrt{ x^{2} + y^{2} } + t^{2} \left( x^{2} + y^{2} \right) $,
\end{enumerate}
and hence
\begin{enumerate}
\item[ ]
$\frac{A}{A+B} \left( a^{2} + b^{2} \right) + \frac{B}{A+B} \left( a^{2} + b^{2} + \sqrt{ \left( a^{2} - c^{2} \right) ^{2} + 4b^{2}(a+c)^{2} } \right) $
\item[ ]
$= \left( a^{2} + b^{2} \right) + \frac{B}{A+B} \sqrt{ \left( a^{2} - c^{2} \right) ^{2} + 4b^{2}(a+c)^{2} }$
\item[ ]
$= t^{2}(t-1)^{2} \left( y^{4} + x^{2}y^{2} \right) + 2t^{2}(t-1)y^{2} \sqrt{ x^{2} + y^{2} } + t^{2} \left( x^{2} + y^{2} \right) $
\item[ ]
$+ \frac{ x^{2} }{ x^{2} + y^{2} } \left( t^{2}(t-1) \left( x^{2} + y^{2} \right) \cdot \left( 2 \sqrt{ x^{2} + y^{2} } + (t-1) \left( x^{2} + y^{2} \right) \right) \right) $
\item[ ]
$= t^{2}(t-1)^{2} \left( x^{2} + y^{2} \right) ^{2} + 2t^{2}(t-1) \left( x^{2} + y^{2} \right) \sqrt{ x^{2} + y^{2} } + t^{2} \left( x^{2} + y^{2} \right) $,
\end{enumerate}
\begin{enumerate}
\item[ ]
$\frac{A}{A+B} \left( b^{2} + c^{2} \right) + \frac{B}{A+B} \left( b^{2} + c^{2} - \sqrt{ \left( a^{2} - c^{2} \right) ^{2} + 4b^{2}(a+c)^{2} } \right) $
\item[ ]
$= \left( b^{2} + c^{2} \right) - \frac{B}{A+B} \sqrt{ \left( a^{2} - c^{2} \right) ^{2} + 4b^{2}(a+c)^{2} }$
\item[ ]
$= t^{2}(t-1)^{2} \left( x^{4} + x^{2}y^{2} \right) + 2t^{2}(t-1)x^{2} \sqrt{ x^{2} + y^{2} } + t^{2} \left( x^{2} + y^{2} \right) $
\item[ ]
$- \frac{ x^{2} }{ x^{2} + y^{2} } \left( t^{2}(t-1) \left( x^{2} + y^{2} \right) \cdot \left( 2 \sqrt{ x^{2} + y^{2} } + (t-1) \left( x^{2} + y^{2} \right) \right) \right) $
\item[ ]
$= t^{2} \left( x^{2} + y^{2} \right) $
\end{enumerate}
and consequently
\begin{enumerate}
\item[ ]
$\left\Vert J \left( f^{(t)} \right) (x,y) \right\Vert ^{2} = \frac{ t^{2}(t-1)^{2} \left( x^{2} + y^{2} \right) ^{2} + 2t^{2}(t-1) \left( x^{2} + y^{2} \right) \sqrt{ x^{2} + y^{2} } + t^{2} \left( x^{2} + y^{2} \right) }{ \left( x^{2} + y^{2} \right) \left( 1 + (t-1) \sqrt{ x^{2} + y^{2} } \right) ^{4} }$
\item[ ]
$= \frac{ t^{2}(t-1)^{2} \left( x^{2} + y^{2} \right) + 2t^{2}(t-1) \sqrt{ x^{2} + y^{2} } + t^{2} }{ \left( 1 + (t-1) \sqrt{ x^{2} + y^{2} } \right) ^{4} }$
\item[ ]
$= \frac{ t^{2} \left( 1 + (t-1) \sqrt{ x^{2} + y^{2} } \right) ^{2} }{ \left( 1 + (t-1) \sqrt{ x^{2} + y^{2} } \right) ^{4} }$
\end{enumerate}
or (equivalently) $$\left\Vert J \left( f^{(t)} \right) (x,y) \right\Vert = \frac{t}{ 1 + (t-1)\sqrt{ x^{2} + y^{2} } }.$$
The same formula is true in case $a<c$ or (equivalently) $x^{2} > y^{2}$. So let us treat the case $a=c$ or (equivalently) $y = \pm x$. Then $$a^{2} + b^{2} \pm 2ab = ( a \pm b )^{2} = \left( t \sqrt{ x^{2} + y^{2} } + t(t-1)y^{2} \pm t(1-t)xy \right) ^{2},$$ which is easily seen to imply that $$\left\Vert J \left( f^{(t)} \right) (x,y) \right\Vert ^{2} = \frac{ \left( t \sqrt{2} \vert x \vert + 2t(t-1)x^{2} \right) ^{2} }{ 2 x^{2} \left( 1 + (t-1) \sqrt{2} \vert x \vert \right) ^{4} }$$ or (equivalently) $$\left\Vert J \left( f^{(t)} \right) (x,y) \right\Vert = \frac{ t }{ 1 + (t-1) \sqrt{2} \vert x \vert }.$$ Finally, since $$J \left( f^{(t)} \right) (0,0) = \left[ \begin{array}{lll} t & 0 \\ \\ 0 & t \end{array} \right] ,$$ it is not difficult to see that $$\left\Vert J \left( f^{(t)} \right) (0,0) \right\Vert = t.$$ We have thus proved that for any $(x,y) \in {\bf R}^{2}$ such that $x^{2} + y^{2} \leq 1$, we have that $$\left\Vert J \left( f^{(t)} \right) (x,y) \right\Vert = \frac{t}{ 1 + (t-1) \sqrt{ x^{2} + y^{2} } }$$ and consequently $$\max\limits_{ x^{2} + y^{2} \leq 1} \left\Vert J \left( f^{(t)} \right) (x,y) \right\Vert = t$$ and $$\min\limits_{ x^{2} + y^{2} \leq 1} \left\Vert J \left( f^{(t)} \right) (x,y) \right\Vert = 1,$$ as it follows by studying the one-variable function $$[0,1] \ni w \mapsto \frac{t}{1+(t-1)w} \in (0, \infty ).$$

\section{The computation of the limit}

For any $t \geq 1$ and for any $z=(x,y) \in {\overline{D}}(0;1)$, we set $$f_{t}(z) = \frac{tz}{ 1 + (t-1) \vert z \vert } = f^{(t)}(x,y),$$ where it is not difficult to see that $f_{1} = id$, while a straightforward computation shows that $$det \left( J \left( f_{t} \right) (z) \right) = \frac{ t^{2} }{ \left( 1 + (t-1) \sqrt{ x^{2} + y^{2} } \right) ^{3} } > 0$$ and it is not difficult to see that on the unit circle we have $f_{t} = id$. The Inverse Function Theorem in 1.53 on page 30 of [2] together with the fact that ${\overline{D}}(0;1)$ is compact implies that $f_{t}$ constitutes a $C^{1}$ diffeomorphism of the closed unit disc, which is the identity on the closed unit circle.
\\ \rm \\
{\bf 4.1. Proposition.} $f_{t} \rightarrow id$ in $\mathcal{G}$ as $t \rightarrow 1^{+}$.
\\ \rm \\
{\bf Proof.} If $t>1$, then a straightforward computation shows that $$\left\Vert f_{t} - id \right\Vert _{ \infty } = \max\limits_{ \vert z \vert \leq 1 } \left\vert f_{t}(z) - z \right\vert = \max\limits_{ \vert z \vert \leq 1 } \frac{ (t-1) \vert z \vert \left( 1 - \vert z \vert \right) }{ 1 + (t-1) \vert z \vert } \leq t-1.$$ So $\left\Vert f_{t} - id \right\Vert _{ \infty } \rightarrow 0$ as $t \rightarrow 1^{+}$. In addition,
\begin{enumerate}
\item[ ]
$\max\limits_{ x^{2} + y^{2} \leq 1 } \left\vert \left( \Re f_{t} \right) _{x}'(x,y) - 1 \right\vert $
\item[ ]
$= \max\limits_{ x^{2} + y^{2} \leq 1 } \left\vert \frac{ t }{ \left( 1 + (t-1) \sqrt{ x^{2} + y^{2} } \right) ^{2} } - 1 + \frac{ t(t-1)y^{2} }{ \sqrt{ x^{2} + y^{2} } \left( 1 + (t-1) \sqrt{ x^{2} + y^{2} } \right) ^{2} } \right\vert $
\item[ ]
$\leq \max\limits_{ x^{2} + y^{2} \leq 1 } \left\vert \frac{ t }{ \left( 1 + (t-1) \sqrt{ x^{2} + y^{2} } \right) ^{2} } - 1 \right\vert $
\item[ ]
$+ t(t-1) \cdot \max\limits_{ x^{2} + y^{2} \leq 1 } \frac{ y^{2} }{ \sqrt{ x^{2} + y^{2} } \left( 1 + (t-1) \sqrt{ x^{2} + y^{2} } \right) ^{2} }$
\item[ ]
$\leq (t-1) + t(t-1) = t^{2} - 1$
\end{enumerate}
and
\begin{enumerate}
\item[ ]
$\max\limits_{ x^{2} + y^{2} \leq 1 } \left\vert \left( \Re f_{t} \right) _{y}'(x,y) - 0 \right\vert $
\item[ ]
$= \max\limits_{ x^{2} + y^{2} \leq 1 } \frac{ t(t-1) \cdot \vert x \vert \cdot \vert y \vert }{ \sqrt{ x^{2} + y^{2} } \left( 1 + (t-1) \sqrt{ x^{2} + y^{2} } \right) ^{2} }$
\item[ ]
$\leq t(t-1) \cdot \max\limits_{ x^{2} + y^{2} \leq 1 } \frac{ \vert x \vert \cdot \vert y \vert }{ \sqrt{ x^{2} + y^{2} } }$
\item[ ]
$\leq t(t-1)$,
\end{enumerate}
since $$\lim\limits_{ (x,y) \rightarrow (0,0) } \frac{ \vert x \vert \cdot \vert y \vert }{ \sqrt{ x^{2} + y^{2} } } = 0,$$
while by analogy $$\max\limits_{ x^{2} + y^{2} \leq 1 } \left\vert \left( \Im f_{t} \right) _{y}'(x,y) - 1 \right\vert \leq t^{2}-1$$ and $$\max\limits_{ x^{2} + y^{2} \leq 1 } \left\vert \left( \Im f_{t} \right) _{x}'(x,y) - 0 \right\vert \leq t(t-1).$$ So $\left\Vert J \left( f_{t} \right) - J(id) \right\Vert _{ \infty } \rightarrow 0$ as $t \rightarrow 1^{+}$ too. \hfill $\bigtriangleup $

\end{document}